\documentclass[a4paper,11pt,DIV=12,
abstract=on
]{scrartcl}
\usepackage{mathtools}
\mathtoolsset{showonlyrefs}
\usepackage{amsmath, amssymb,amsthm,xargs}
\usepackage{fontenc}
\usepackage[utf8]{inputenc}
\usepackage[english]{babel}
\usepackage{graphicx}
\usepackage{subcaption}
\usepackage{algorithm, booktabs}
\usepackage[noend]{algpseudocode}
\usepackage{enumitem,listings,xcolor}

\newtheorem{theorem}{Theorem}[section]

\setkomafont{sectioning}{\rmfamily\bfseries}
\setkomafont{title}{\rmfamily}
\usepackage{hyperref}

\DeclareCaptionLabelSeparator{periodspace}{.\ }
\newcommand{\R}{\mathbb{R}}

\newcommand{\Z}{\mathbb{Z}}

\newcommand{\tT}{\mathrm{T}}
\newcommandx{\abs}[2][1=\@empty]{#1\lvert #2 #1\rvert}
\newcommandx{\norm}[3][1=\@empty,3=2]{#1\lVert #2 #1\rVert_{#3}}
\DeclareMathOperator*{\argmin}{arg\,min}
\DeclareMathOperator{\prox}{prox}

\DeclareMathOperator{\sgn}{sgn}

\hypersetup{
pdfauthor={Ronny Bergmann, Andreas Weinmann},
pdftitle={Inpainting of Cyclic Data using First and Second Order Differences},
pdfsubject={Submitted Conference Proceedings Paper},%
pdfcreator = {pdflatex and TextMate},%
pdfkeywords={Inpainting,%
	variational models with higher order derivatives,%
	cyclic data, phase-valued data,%
	cyclic proximal point algorithm},
plainpages=false, pdfstartview=FitH, pdfview=FitH, pdfpagemode=UseOutlines,%
bookmarksnumbered=true,bookmarksopen=false,bookmarksopenlevel=0,%
colorlinks=true,linkcolor=black,citecolor=black,urlcolor=black%
}
\title{Inpainting of Cyclic Data using \\
First and Second Order Differences}
\author{%
Ronny Bergmann\thanks{Department of Mathematics,
Technische Universität Kaiserslautern,
Paul-Ehrlich-Str. 31, 67663 Kaiserslautern, Germany,
bergmann@mathematik.uni-kl.de.}
, 
Andreas Weinmann\thanks{Department of Mathematics,
Technische Universität München and
Fast Algorithms for Biomedical Imaging Group, Helmholtz-Zentrum München,
Ingolstädter Landstr. 1, 85764 Neuherberg, Germany,
andreas.weinmann@tum.de.}
}
\date{\today}
\begin{document}
\maketitle
\begin{abstract}
	\parindent0mm
	\parskip0mm
	\noindent
Cyclic data arise in various image and signal processing applications such as
interferometric synthetic aperture radar, electroencephalogram data analysis,
and color image restoration in HSV or LCh spaces.
In this paper we introduce a variational inpainting model for cyclic data which
utilizes our definition of absolute cyclic second order differences.
Based on analytical expressions for the proximal mappings of these differences
we propose a cyclic proximal point algorithm (CPPA) for minimizing the
corresponding functional.
We choose appropriate cycles to implement this algorithm in an efficient way.
We further introduce a simple strategy to initialize the unknown inpainting
region. Numerical results both for synthetic and real-world data  demonstrate
the performance of our algorithm.
\end{abstract}

\paragraph{Keywords.}\noindent
Inpainting,
variational models with higher order differences,
cyclic data,
phase-valued data,
cyclic proximal point algorithm.
%
\section{Introduction}\label{sec:Introduction}
Image inpainting is a frequently arising problem in image processing.
Examples are restoring scratches in photographs, removal of superimposed
objects, dealing with areas removed by a user, digital zooming, edge decoding,
restoration of defects in audio/video recordings or in seismic data.
The term `inpainting' first appeared in~\cite{BSCB00}, but earlier work on
disocclusions was already done,  e.g.,
in~\cite{caselles1998axiomatic,masnou1998level}.
In this respect also interpolation, approximation, and extrapolation problems
may be viewed as inpainting problems.
Inpainting is a very active field of research which has been tackled by various
approaches.
For a good overview we refer to the (tutorial)
papers~\cite{BBCS10,CDOS12,chan2005image,guillemot2014image}.
While exemplar-based and sparsity-based (dictionary/frame/tensor) methods
are in general better suited for filling large texture areas,
diffusion-based and corresponding variational techniques show good results for
natural images.
The total variation (TV) regularized model proposed in~\cite{ROF92} for
denoising was first applied to inpainting in~\cite{BBCSV01,CS01}.
It was later also used in combination with other methods, however, the TV
regularizer typically introduces a staircasing effect in the
corresponding minimizer.
A simple method to avoid these artifacts consists in the incorporation of
second order derivatives into the model.
Indeed, starting with~\cite{CL97} various approaches with higher order
derivatives have been proposed, see,
e.g.,~\cite{BKP09,CMM00,DWB09,HS06,LBU2012,LLT03,PS13,Sche98,SS08,SST11,VBK13}.
In this paper, we address the problem of inpainting cyclic data using a 
variational model with second order cyclic differences.
In general, manifold-valued data processing has recently gained a lot of
interest.
Examples are wavelet-type multiscale transforms for manifold
data~\cite{GW09,RDSDS05,Wein12} and manifold-valued partial differential
equations~\cite{CTDF04,GHS13}. Also statistical issues on
Riemannian manifolds have been considered~\cite{Fle13,FJ07,Pen06},
in particular the statistics of circular data~\cite{fisher95,JS2001}.

\paragraph{Related work.} \noindent%
Although very popular for processing images
with scalar and vector-valued data, TV minimization has only very recently been
applied to cyclic structures.
From a theoretical point of view TV functionals for manifold-valued functions
have been studied in~\cite{GM06,GM07}.
These papers extend the previous work~\cite{GMS93} on $\mathbb S^1$-valued
functions where, in particular, the existence of minimizers of certain energies
is shown in the space of functions with bounded total cyclic variation.
First order TV minimization for cyclic data in image processing has been
investigated in~\cite{SC11,CS13}.
The authors unwrap the data to the real line and propose an algorithm based on
functional lifting which takes the periodicity into account.
In particular, they also consider cyclic inpainting.
An algorithm for TV minimization on Riemannian manifolds was proposed
in~\cite{LSKC13}.
The approach is based on a reformulation as a multilabel optimization problem
with an infinite number of labels.
Using convex relaxation techniques, the resulting hard optimization problem is
approximated which also requires the discretization of the manifold.
Another approach for denoising manifold-valued data via first order TV
minimization was given in~\cite{WDS2013}. The authors propose cyclic and
parallel proximal point algorithms which will also be our method of choice.

\paragraph{Contributions.} \noindent%
We propose two models for inpainting of cyclic data using first and second
order absolute cyclic differences.
In our preprint~\cite{BLSW14} we introduced absolute second
order differences for cyclic data in a sound way.
We further deduced analytical expressions for the proximal mappings of these
differences.
Here, our first model considers the noise free inpainting situation, whereas the
second one handles simultaneously inpainting and denoising.
The variational formulations allow for the decomposition of the whole
functionals into simpler ones, for each of which the proximal mappings are
given explicitly.
Thus, the minimizers can be computed efficiently by a cyclic proximal point
method.
We propose a suitable initialization of the inpainting area.
We demonstrate by numerical examples the strength of our algorithm.
Compared to~\cite{SC11,CS13} we neither have to employ Fréchet means nor to
discretize the manifold.

\paragraph{Organization.}
In Sec.~\ref{sec:cyc_diff} we introduce our absolute second order cyclic
differences and provide analytical expressions for their proximal mappings.
Then, in Sec.~\ref{sec:model}, we introduce our inpainting model
and propose a procedure to initialize the unknown inpainting region.
Sec.~\ref{sec:cpp} describes the cyclic proximal point algorithm.
Finally, Sec.~\ref{sec:numerics} contains numerical examples.
Conclusions and directions of future work are given in
Sec.~\ref{sec:conclusions}.
%
%
\section{Absolute First and Second Order Cyclic Differences and Their Proximal Mappings} \label{sec:cyc_diff}
\begin{figure}[tbp]\centering
	\includegraphics{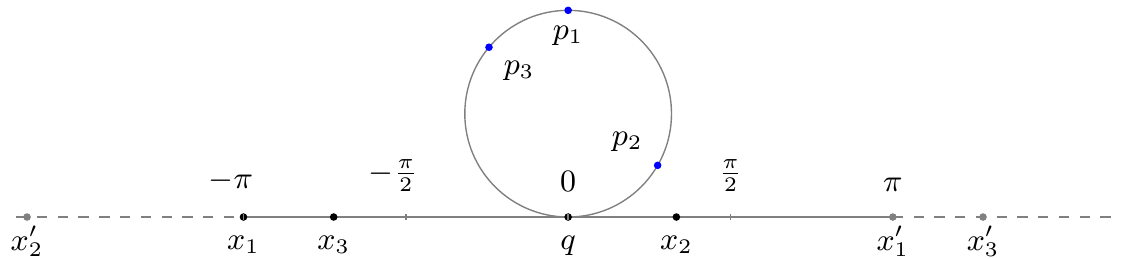}
 \caption{Three points  \(p_j\), $j=1,2,3$, on the circle and their possible unwrappings \(x_1,x_2,x_3 \in [-\pi,\pi)\) with respect to the origin \(q\) and
 other possibilities \(x'_2,x_1,x_3\) and \(x_2,x_3,x'_1\) that correspond to the same situation on \(\mathbb S^1\). These are taken into account for \(\text{d}(x;w)\).}\label{fig:circle}
\end{figure}

%
Let \(\mathbb S^1:= \{p_1^2 + p_2^2 = 1: p = (p_1,p_2)^\tT \in \R^2\}\) be the
unit circle endowed with the {\itshape geodesic distance}
\( d_{\mathbb S^1} (p,q) := \arccos( \langle p,q \rangle )\).
Given a base point \(q \in \mathbb S^1\), the {\itshape exponential map}~
\(\exp_q: \R \rightarrow \mathbb S^1\) from the tangent
space~\(T_q\mathbb S^1 \simeq \R\) of \(\mathbb S^1\) at \(q\) onto
\(\mathbb S^1\) is defined by
\[
	\exp_q(x) = R_x q, \qquad R_x :=
		\begin{pmatrix}
			\cos x & -\sin x\\
			\sin x & \cos x
		\end{pmatrix}.
\]
This map is $2\pi$-periodic, i.e., $\exp_q(x) = \exp_q((x)_{2\pi})$ for any
$x \in \R$, where $(x)_{2\pi}$ denotes the unique point in $[-\pi,\pi)$ such
that $x = 2\pi k + (x)_{2\pi}$, $k \in \mathbb Z$.
For~$p,q \in \mathbb S^1$ with $\exp_q(0) = q$, there is a unique
$x \in [-\pi,\pi)$ satisfying~$\exp_q(x) = p$.
Given such representants $x_j \in [-\pi,\pi)$ of
$p_j \in \mathbb S^1$, $ j = 1,2$ centered at an arbitrary base
point~$q\in \mathbb S^1$ the geodesic distance becomes
\[
	d_{\mathbb S^1}(p_1,p_2) = d(x_1,x_2) = \min_{k \in \Z }\lvert x_2 - x_1 + 2\pi k\rvert
	= \lvert(x_2 - x_1)_{2\pi}\rvert 
\]
which is of course independent of $q$.
We want to define higher order differences for points
$(p_j)_{j=1}^d \in (\mathbb S^1)^d$ using their representants
$x := (x_j)_{j=1}^d \in [-\pi,\pi)^d$.
To achieve independence of the base point the differences must be shift
invariant modulo $2\pi$, see Fig.~\ref{fig:circle}.
Let $1_d$ denote the vector with $d$ entries 1.
We define the {\itshape absolute cyclic difference} of $x \in [-\pi,\pi)^d$ 
with respect to a difference filter
$w \in \mathbb R^d$ with $\langle w, 1_d \rangle = 0$ by
\begin{equation} \label{def_cyc_diff}
	\mathrm{d}(x;w) 
	:= \min_{\alpha \in \R} \langle [x + \alpha 1_d]_{2 \pi}, w  \rangle,
\end{equation}
where $[x]_{2\pi}$ denotes the componentwise application of
$(t)_{2\pi}$ if $t \not = (2k +1) \pi $, $k \in \Z$ and 
$[(2k +1) \pi]_{2\pi} = \pm \pi$, $k \in \Z$.
Let 
\[
b_1 := (-1,1)^\tT \quad \text{and} \quad  b_2 := (1,-2,1)^\tT, \; b_{1,1} := (-1,1,1,-1)^\tT
\]
be a first order (forward) difference filter, and two second order difference
filters, respectively.
For $w \in {\cal B} := \{b_1,b_2,b_{1,1} \}$ we have shown in our accompanying
preprint~\cite{BLSW14} that the absolute cyclic differences can be rewritten as
\begin{equation} \label{sec_diff}
	\mathrm{d}(x;w) 
		= \left( \langle x,w \rangle \right)_{2\pi}.
\end{equation}
Clearly, we have $\mathrm{d}(x;b_1) = d(x_1,x_2)$.
Interestingly, the definition \eqref{def_cyc_diff} and \eqref{sec_diff} do not
coincide, e.g., for third order cyclic differences~\cite{BLSW14}.

Next we are interested in proximal mappings of absolute cyclic
differences. Recall that for a proper, closed, convex function 
$\varphi: \R^N \rightarrow (-\infty,+\infty]$
and $\lambda > 0$ the {\itshape proximal mapping} $\prox_{\lambda \varphi}: \R^N \rightarrow \R^N$
is well defined by
\begin{equation}\label{prox_R}
	\prox_{\lambda \varphi}(f)
	:=
	\argmin_{x \in \R^N} \frac{1}{2} \lVert f-x\rVert_2^2 + \lambda \varphi (x).
\end{equation}
We introduce the proximal mapping 
$
\prox_{\lambda \mathrm{d}(\cdot;w)}:(\mathbb S^1)^d \rightarrow (\mathbb S^1)^d
$ 
by
\[
\prox_{\lambda \mathrm{d}(\cdot;w)}(f)
:= \argmin_{x \in [-\pi,\pi)^d } \frac{1}{2} \sum_{j=1}^d d(x_j,f_j) ^2
	+ \lambda \mathrm{d}(x;w), \quad \lambda > 0.
\]
The following theorem determines the proximal mapping analytically
for $w \in  {\cal B}$.
In particular, the mapping is single-valued for $f \in [-\pi,\pi)^d$ with
$\lvert(\langle f,w \rangle)_{2\pi} \rvert < \pi$ and two-valued for
$\lvert(\langle f,w \rangle)_{2\pi} \rvert = \pi$.
Note that for $w = b_1$ the second case appears exactly if $f_1$ and
$f_2$ are antipodal points.
For a proof we refer to our preprint~\cite{BLSW14}.
%
\begin{theorem} \label{lem:proxy_b1}
For $w \in {\cal B}$ set $s := \sgn(\langle f,w \rangle)_{2 \pi}$.
Let $f \in [-\pi,\pi)^d$, where $d$ is adapted to the respective
length of $w$, $\lambda > 0$, and $m:= \min \left\{\lambda,%
\frac{\lvert (\langle f,w \rangle)_{2\pi} \rvert}%
{\lVert w\rVert_2^2} \right\}$.
\begin{itemize}
\item[i)]
If $\lvert(\langle f,w \rangle)_{2\pi} \rvert < \pi$, then 
$
\prox_{\lambda \mathrm{d}(\cdot;w)}(f) =  (f - s \, m \,w)_{2\pi}$.

\item[ii)]
If $\lvert(\langle f,w \rangle)_{2\pi} \rvert = \pi$,
then 
$
\prox_{\lambda \mathrm{d}(\cdot;w)}(f) 
= \{(f + s \, m \, w)_{2\pi}, (f - s \, m \, w)_{2\pi}\}.
$
\end{itemize}
\end{theorem}
%
For handling noisy data we will further need the following proximal mapping:
%
\begin{theorem} \label{Theo:Prox_quad}
For $f,g \in [-\pi,\pi)^N$ we have
\begin{align}
	\prox_{\lambda d(\cdot,f)}(g) &:= \argmin_{x}
	\sum_{j=1}^{N} \left( d(g_j,x_j)^{2}+\lambda d(f_j,x_j)^{2} \right)\\
&=
\left( \frac{g+\lambda f}{1+\lambda} + \frac{\lambda}{1+\lambda} \, 2\pi \, v \right)_{2\pi},
\end{align}
where \(d(g,f) := \displaystyle\sum_{j=1}^N d(g_j,f_j)\) and $v = (v_j)_{j=1} ^N \in \R^N$ is defined by
\[
v_j :=
\begin{cases} 
	0 & \mbox{ if \(\lvert g_j-f_j\rvert \le \pi \)},\\
	\sgn(g_j - f_j) & \mbox{ if \(\lvert g_j-f_j\rvert > \pi\)}.
\end{cases}
\]
\end{theorem}

\section{Inpainting Models for Cyclic data} \label{sec:model}
%
Given an image domain \(\Omega_0 = \{1,\ldots,N\} \times \{1,\ldots,M\}\),
the inpainting region $\Omega\subset\Omega_0$ is the subset where the pixel
values \(f_{i,j}\), \((i,j) \in \Omega\) are unknown.
The (noiseless) inpainting problem consists of finding a function $x$ on
\(\Omega_0\) from data $f$ given on $\bar \Omega = \Omega_0 \backslash \Omega$
such that $x$ is a suitable extension of~$f$ to $\Omega_0$. Let 
$d_2 (x) := \mathrm{d}(x;b_2)$ and $d_{1,1} (x) := \mathrm{d}(x;b_{1,1})$.
Our functional for inpainting of noiseless cyclic data reads
\begin{equation}\label{eq:2DTVfunctionalInpainting}
	\begin{split}
		\argmin_{{x}\in [-\pi,\pi)^{N,M}}
			\alpha \operatorname{TV}_{1}^{\Omega}(x)
			+ \beta \operatorname{TV}_{2}^{\Omega}(x)
			+ \gamma \operatorname{TV}_{1,1}^{\Omega}(x),\\
			\text{ s.t. }\quad x_{i,j}=f_{i,j}
			\quad\text{ for all }\quad (i,j)\in \bar\Omega,
	\end{split}
\end{equation}
where $\alpha := (\alpha_1,\alpha_2,\alpha_2,\alpha_4)$, $\beta := (\beta_1,\beta_2)$ and
the restricted first and second order difference terms given by 
\begin{align}
	\alpha\operatorname{TV}_{1}^{\Omega}(x)
	&=
		\alpha_1 \sum_{(i,j)} d(x_{i,j},x_{i+1,j})  
		+
		 \alpha_2 \sum_{(i,j)\}} d(x_{i,j},x_{i,j+1})
\\
	&\quad
		+
		\frac{1}{\sqrt{2}} \bigg( \alpha_3\sum_{(i,j)} d(x_{i,j},x_{i+1,j+1})
		+
		\alpha_4 \sum_{(i,j)} d(x_{i,j+1},x_{i+1,j}) \bigg),
	\label{eq:2DTVInpainting}
\end{align}
\begin{align}
	\beta\operatorname{TV}_{2}^{\Omega}(x)
	&=
	\beta_1 \sum_{(i,j)}
	d_2(x_{i-1,j},x_{i,j},x_{i+1,j})
		+
	\beta_2\sum_{(i,j)}
	d_2(x_{i,j-1},x_{i,j},x_{i,j+1}),
	\label{eq:TV2isoInpainting}
\end{align}
and
\begin{align}\label{eq:TV2mixInpainting}
	\gamma \operatorname{TV}_{1,1}^{\Omega}(x)
	&=\gamma
	\sum_{(i,j)}
	d_{1,1}(x_{i,j},x_{i+1,j},x_{i,j+1},x_{i+1,j+1}),		
\end{align}
where the sums are taken only for those $(i,j)$ for which at last one entry $x_{a,b}$ in the corresponding differences is contained in $\Omega$. We use the notation \(\operatorname{TV}\) since the model of the first order differences resembles an anisotropic TV model.

For the inpainting problem in the presence of noise the requirement of equality on $\bar \Omega$ is replaced by $x$ being an approximation of $f$: 
\begin{align}\label{eq:2DTVfunctionalInpaintingWithNoise}
\argmin_{{x}\in [-\pi,\pi)^{N,M}}
	F_{\bar \Omega}(x; f)
	+ \alpha \operatorname{TV}_{1}(x)
	+ \beta \operatorname{TV}_{2}(x)
	+ \gamma \operatorname{TV}_{1,1}(x),
\end{align}
where
\[
F_{\bar \Omega}(x; f) := \sum_{(i,j)\in \bar \Omega} d(x_{i,j},f_{i,j})^2.
\]
and the first and second order difference terms sum over all indices in $\Omega_0$ now.

\paragraph{Initialization of the inpainting region.}
Since the inpainting problem does not posses a unique minimizer the
initialization of the inpainting area is crucial.
We present a method which is related to the idea of unknown boundary conditions
used by Almeida and Figueiredo in~\cite{AF13}.
It can also be viewed as an implicit version of the ordering method of pixels
by adapted distance functions used by März
in~\cite{marz2011image,marz2013well}.
To this end, we initialize \(x_{i,j} = f_{i,j}\) for \((i,j)\in\bar \Omega\). The other ones are considered as not initialized.
We use first, second and mixed order differences \(d=d_1, d_2\) and \(d_{1,1}\) and let \(t\in\{1,2,(1,1)\}\). Let \(x := (x_{k_1},\ldots,x_{k_l})^\tT\) be a set of points corresponding to a stencil of such a difference term~\(d_t\). If \(k_i\in\Omega\), \(i\in\{1,\ldots,l\}\), is the unique index such that \(x_{k_i}\) is not yet initialized, i.e.,
there is exactly one unknown point at \(k_i\in\Omega\) in \(x\), we can
initialize this value as follows. The minimal value for the absolute cyclic
finite difference is 
\(
	0 = \left( \langle x, b_t \rangle \right)_{2\pi}
\)
and this equation provides an initial value for \(x_{k_i}\). Such a situation
of exactly one unknown index \(k_i\) always exists at the boundary of the
initialized area.
%
\section{Cyclic Proximal Point Algorithm} \label{sec:cpp}
%
Since the proximal mappings of our absolute cyclic differences can be
efficiently computed using their analytical expressions
in Theorem \ref{lem:proxy_b1} and Theorem \ref{Theo:Prox_quad}, we suggest to
apply a cyclic proximal point algorithm to find a minimizer for the
inpainting problem. 
Recently, the proximal point algorithm (PPA) on the Euclidean
space~\cite{Roc76} was extended to Riemannian manifolds of non-positive
sectional curvature~\cite{FO02} and also to Hadamard spaces~\cite{Bac13}.
A cyclic PPA (CPPA) on the Euclidean space was given in~\cite{Ber10,Ber11}
and on Hadamard spaces in~\cite{Bac13a}.
Unfortunately, one of the simplest manifolds that is not of Hadamard type is
the circle $\mathbb S^1$.
However, under certain assumptions we were able to prove the convergence of
the CPPA to a minimizer of the denoising problem for cyclic data, see~\cite{BLSW14}.
A similar proof can also be given for the inpainting problem. 
Indeed, we have observed convergence of our algorithm in all numerical tests. 

In the CPPA the original function $J$ is split into a sum
$J = \sum_{l=1}^c J_l$ and the proximal mappings of the functions $J_l$ are
applied in each iteration cycle, i.e.,
\[
x^{(k+1)} = \prox_{\lambda_kJ_c}
	\left( \prox_{\lambda_k J_{c-1}} 
		\left( \ldots \left( \prox_{\lambda_k J_{1}} (x^{(k)}) \right)\right)
	\right).
\]
For $J = J_1 + J_2$, where $J_1,J_2: \mathbb R^N \rightarrow (-\infty,+\infty]$ 
are proper, closed convex functions, it is well known that  the nested PPA 
\[
x^{(k+1)} = \prox_{\lambda J_2} \left( \prox_{\lambda J_{1}}(x^{(k)}) \right)
\]
converges for any fixed parameter $\lambda >0$ to a fixed point of
$\prox_{\lambda J_2} \circ \prox_{\lambda J_{1}}$.
Unfortunately this fixed point is  not a minimizer of $J$ but of
$J_2 + {}^\lambda \! J_{1}$, where ${}^\lambda \! J_{1}$ denotes the Moreau
envelope of $J_1$.
Convergence to the correct minimizer can be achieved by choosing an iteration
dependent sequence  $\{\lambda_k\}_k$ fulfilling
 \begin{equation}\label{eq:CPPAlambda}
	\sum_{k=0}^\infty \lambda_k = \infty,
		\quad \text{and}
		\quad \sum_{k=0}^\infty \lambda_k^2 < \infty,
\end{equation}
see~\cite{Bac13a,Ber11}. A specific splitting of our inpainting model \eqref{eq:2DTVfunctionalInpainting} for the CPPA is given in the appendix.
\section{Numerical Results} \label{sec:numerics}
\begin{figure}[tbp]\centering
	\begin{subfigure}[t]{.31\textwidth}\centering
		\includegraphics{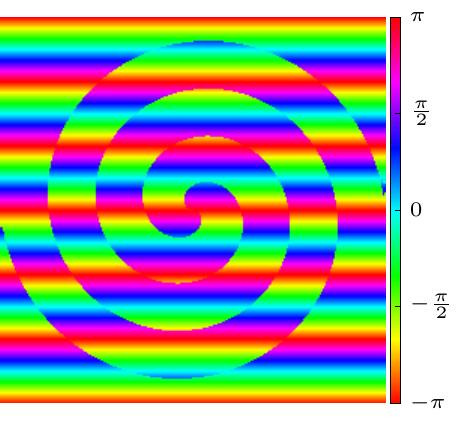}
		\caption{Original image,\\\(257\times 257\) pixel.}\label{subfig:sar-orig}
	\end{subfigure}
	\begin{subfigure}[t]{.31\textwidth}\centering
		\includegraphics{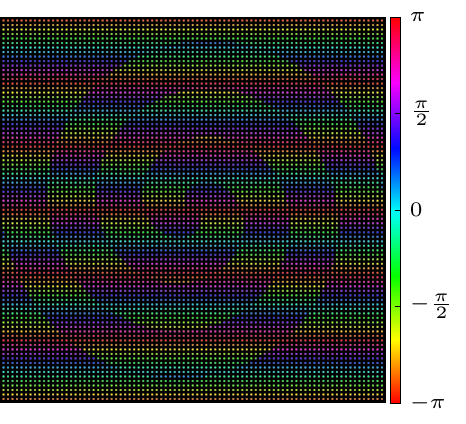}
		\caption{Masked image,\\89.8\% data lost.}\label{subfig:sar-mask}
	\end{subfigure}
	\begin{subfigure}[t]{.31\textwidth}\centering
		\includegraphics{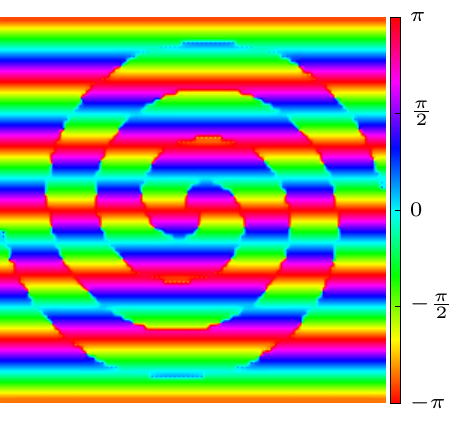}
		\caption{Inpainting result,\\\(\alpha = (1,1,0,0)^\tT\),\\\(\beta = (1,1)^\tT\), \(\gamma=1\).}\label{subfig:sar-result}
	\end{subfigure}
	\caption{The synthetic SAR data in~\subref{subfig:sar-orig} is reduced by a factor of nine by removing two thirds of all rows and columns. They are indicated as black pixels in~\subref{subfig:sar-mask}, the right
	image~\subref{subfig:sar-result} shows the reconstruction based on first and second order cyclic differences.}
	\label{fig:sar}
\end{figure}

\begin{figure}[tbp]\centering
	\begin{subfigure}[t]{.31\textwidth}\centering
		\includegraphics{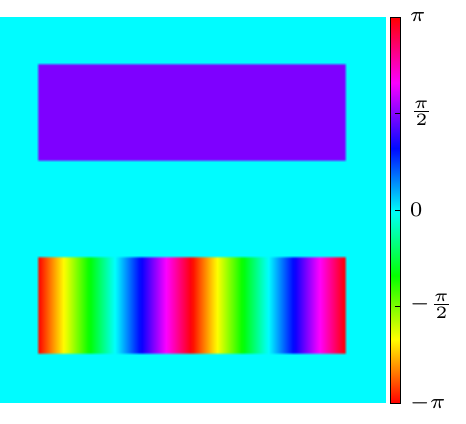}
		\caption{Original image,\\\(128\times 128\) pixel.}\label{subfig:blocks-orig}
	\end{subfigure}
	\begin{subfigure}[t]{.31\textwidth}\centering
		\includegraphics{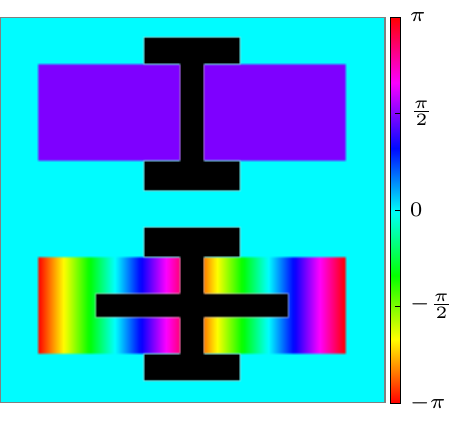}
		\caption{Masked image.}\label{subfig:blocks-mask}
	\end{subfigure}
	\begin{subfigure}[t]{.31\textwidth}\centering
		\includegraphics{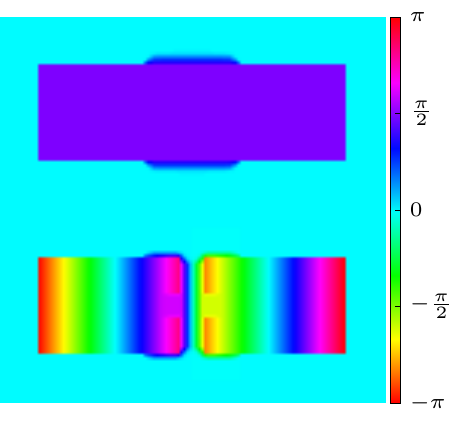}
		\caption{\(\mathbb R\)-valued inpainting\\with parameters of~\subref{subfig:blocks-TV12}.}\label{subfig:blocks-R}
	\end{subfigure}
	\\
	\begin{subfigure}[t]{.31\textwidth}\centering
		\includegraphics{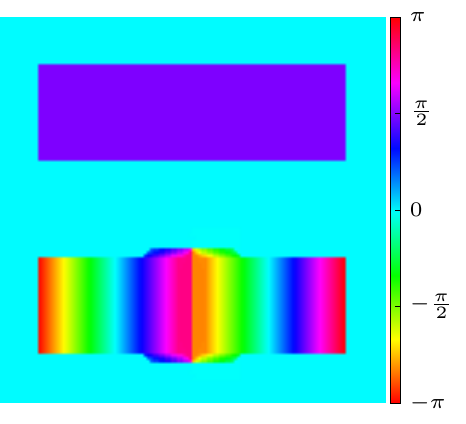}
		\caption{\(\mathbb S^1\)-valued inpainting
		\\\(\alpha=(2,2,2,2)^\tT\).}\label{subfig:blocks-TV}
	\end{subfigure}
	\begin{subfigure}[t]{.31\textwidth}\centering
			\includegraphics{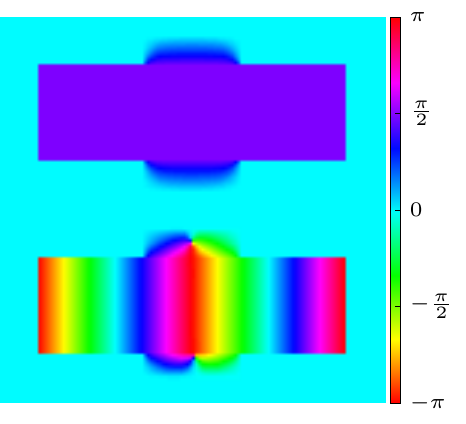}
		\caption{\(\mathbb S^1\)-valued inpainting
		\\ \(\beta=(1,1)^\tT\), \(\gamma=1\).}\label{subfig:blocks-TV2}
	\end{subfigure}
	\begin{subfigure}[t]{.31\textwidth}\centering
		\includegraphics{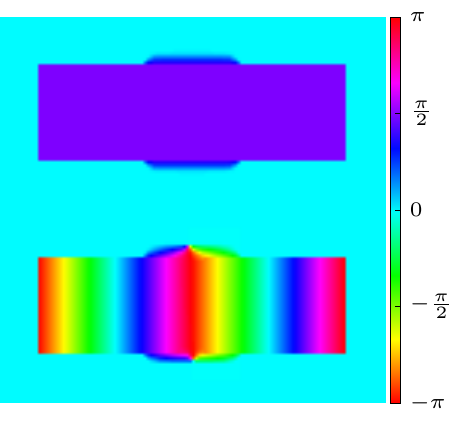}
		\caption{\(\mathbb S^1\)-valued inpainting,\\\(\alpha=(2,2,2,2)^\tT\),\\ \(\beta=(1,1)^\tT\), \(\gamma=1\).}\label{subfig:blocks-TV12}
	\end{subfigure}
	\caption{
	%
Inpainting with first and second order differences: \subref{subfig:blocks-orig} from the original image \subref{subfig:blocks-mask} some parts (black) are lost. \subref{subfig:blocks-R} A real-valued inpainting fails; \subref{subfig:blocks-TV} a first order model reconstructs the constant region perfectly; \subref{subfig:blocks-TV2} a pure second order model has linear artifacts; \subref{subfig:blocks-TV12} a first and second order model performs best.}\label{fig:block}
\end{figure}

For the numerical computations of the following examples, the presented algorithms were implemented in \textsc{Matlab}. The experiments were performed
on a MacBook Pro with an Intel Quad Core i5, 2.6\,Ghz and 8\,GB of RAM
on OS X 10.9.2.
\paragraph{Interpolation and Approximation.}
As a first example we consider an synthetic SAR data sample taken
from~\cite{GP98}\footnote{online available at
\url{ftp://ftp.wiley.com/public/sci_tech_med/phase_unwrapping/data.zip.}}, see
Fig.~\ref{fig:sar}\,(\subref{subfig:sar-orig}). We destroy about 89.8\% of the
data by removing all but the rows and columns that are not divisible by 3, see
Fig.~\ref{fig:sar}\,(\subref{subfig:sar-mask}). This is taken as input for the
CPPA in order to minimize~\eqref{eq:2DTVfunctionalInpainting} using the
parameters \(\alpha=(1,1,0,0)^\tT\), \(\beta = (1,1)^\tT\) and \(\gamma = 1\).
The result is shown in Fig.~\ref{fig:sar}\,(\subref{subfig:sar-result}), where
the linear parts are reconstructed perfectly, while the edges are interpolated
and hence suffer from linearization of the original circular edge path.
The runtime is about 80 seconds for the image of size \(257\times257\) pixel
when using \(k=700\) iterations as a stopping criterion for the CPPA from
Sec.~\ref{sec:cpp}.
%
\paragraph{Inpainting for Restoring Image Regions.}
A main application of inpainting is to restore destroyed image regions in
noiseless images. We use the first model~\eqref{eq:2DTVfunctionalInpainting}
and consider an example adapted from~\cite{PS13}, where a similar image was
used to demonstrate regularization with a second order model for real valued
images.
We extend their experiment by including a region with linear increase that is
wrapped twice, cf. Fig.~\ref{fig:block}\,(\subref{subfig:blocks-orig}).
We remove a vertical strip in the middle of both regions and stripes between
the fore- and background. Furthermore for the second, linearly increasing
region we mask a small band in the middle, cf.
Fig.~\ref{fig:block}\,(\subref{subfig:blocks-mask}).

We then employ a real valued inpainting using first and second order
differences, cf. Fig.~\ref{fig:block}\,(\subref{subfig:blocks-R}). The constant
rectangle shows a similar behavior to~\cite{PS13}, where the smoothing at the
top and bottom is reduced here. This is due to employment of both first and
second order real valued differences. Most noticeably, the linearly increasing
region is not reconstructed.

The
Figs.~\ref{fig:block}\,(\subref{subfig:blocks-TV})--(\subref{subfig:blocks-TV12})
illustrate the effects of first and second order absolute cyclic differences.
Fig.~\ref{fig:block}\,(\subref{subfig:blocks-TV}) uses only the first order
model, (\subref{subfig:blocks-TV2}) only the second order differences. and
(\subref{subfig:blocks-TV12}) combines both.
The first order absolute cyclic differences reconstruct the constant region
perfectly, but also produce the well known staircasing in the lower part.
The second order cyclic model introduces a smooth transition between fore- and
background. However, it perfectly reconstructs the linear increase. Combining
both the first and second order cyclic models yields a perfect reconstruction
of the linearly increasing region while reducing the smooth transition,
cf. Fig.~\ref{fig:block}\,(\subref{subfig:blocks-TV12}).

\begin{figure}[tbp]\centering
	\begin{subfigure}[t]{.31\textwidth}\centering
		\includegraphics{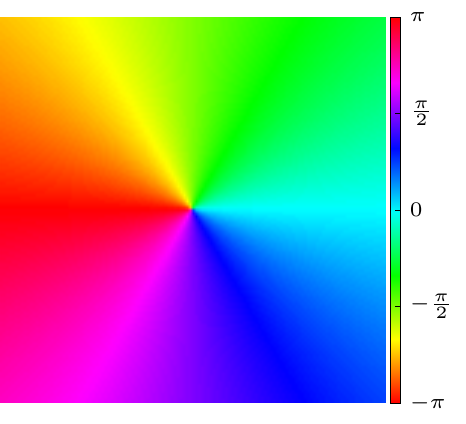}
			\caption{Original image,\\\(128\times 128\) pixel.}\label{subfig:cont:orig}
		\end{subfigure}
		\begin{subfigure}[t]{.31\textwidth}\centering
			\includegraphics{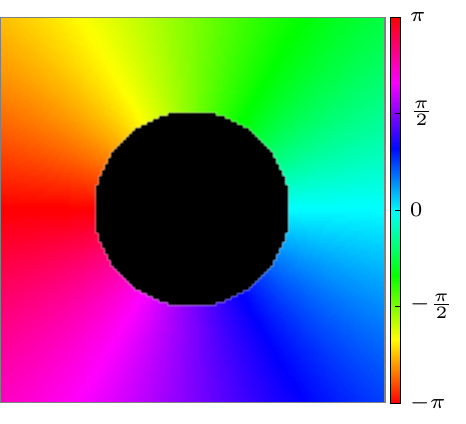}
			\caption{Masked image.
			}\label{subfig:cont:mask}
		\end{subfigure}
		\\
		\begin{subfigure}[t]{.31\textwidth}\centering
		\includegraphics{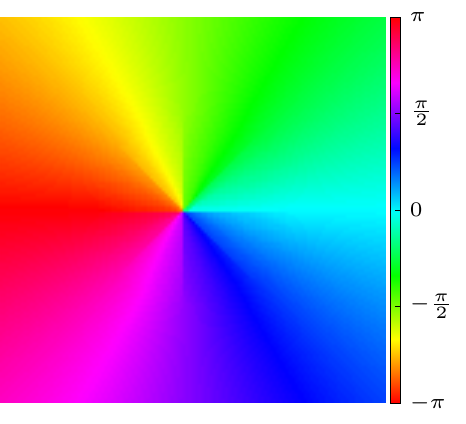}
			\caption{Inpainting,\\\(\alpha=\tfrac{1}{2}(1,1,1,1)^\tT\),\\\(\beta = (0,0)^\tT\), \(\gamma=0\).}\label{subfig:cont:inp1}
		\end{subfigure}
		\begin{subfigure}[t]{.31\textwidth}\centering
			\includegraphics{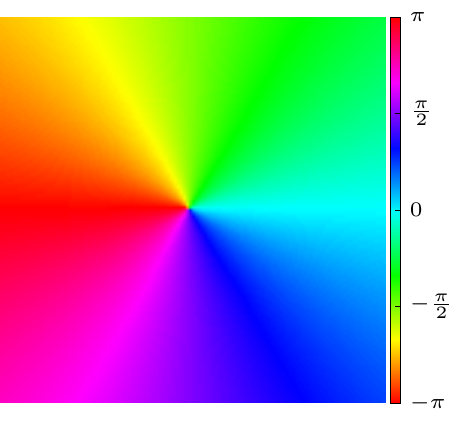}
			\caption{Inpainting,\\ \(\alpha = \tfrac{1}{2}(1,1,1,1)^\tT\),\\\(\beta = \tfrac{1}{4}(1,1)^\tT\), \(\gamma=\tfrac{1}{4}\).}\label{subfig:cont:inp2}
		\end{subfigure}
		\caption{We mask a circular region at the center
		of~\subref{subfig:cont:orig}, see~\subref{subfig:cont:mask}.
		The reconstruction used in~\subref{subfig:cont:inp1} employs only first
		order cyclic differences and produces staircasing. Combining first and
		second order cyclic differences in~\subref{subfig:cont:inp2} we obtain
		a nearly perfect reconstruction.}\label{fig:cont}
\end{figure}

As a second reconstruction example we consider the function
\(\operatorname{atan2}(y,x)\) sampled on a regular grid in
\([-\tfrac{1}{2},\tfrac{1}{2}]^2\) having \(128\) sampling points in each
dimension, cf. Fig.~\ref{fig:cont}\,(\subref{subfig:cont:orig}). We take a
circular mask in the center of the image, see
Figure~\ref{fig:cont}\,(\subref{subfig:cont:mask}), where the mask is shown in
black. In this experiment we compare the results using only first order
absolute cyclic differences with a combined approach of first and second order
cyclic model.

When only using first order differences, we obtain a result that again reveals
staircasing, cf. Fig.~\ref{fig:cont}\,(\subref{subfig:cont:inp1}). It prefers
\(x\)- and \(y\)-axis, and both diagonals, which can be seen by the crosses
created in the middle.
By also including second order differences we obtain almost the original image, cf. Fig.~\ref{fig:cont}\,(\subref{subfig:cont:inp2}).

\begin{figure}[tbp]\centering
		\begin{subfigure}[t]{.49\textwidth}\centering
			\includegraphics{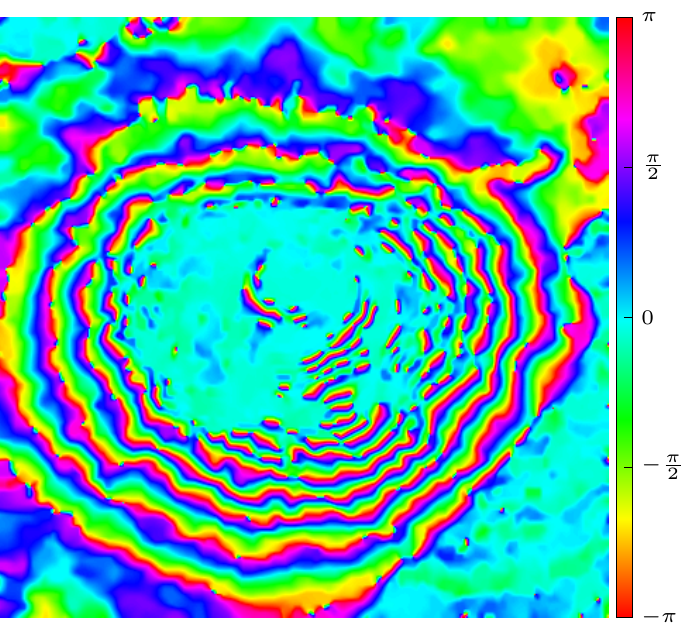}
			\caption{Denoised Mt. Vesuvius image,\\\(432\times 426\) pixels.}\label{subfig:Vesuv:Denoised}
		\end{subfigure}
		\begin{subfigure}[t]{.49\textwidth}\centering
			\includegraphics{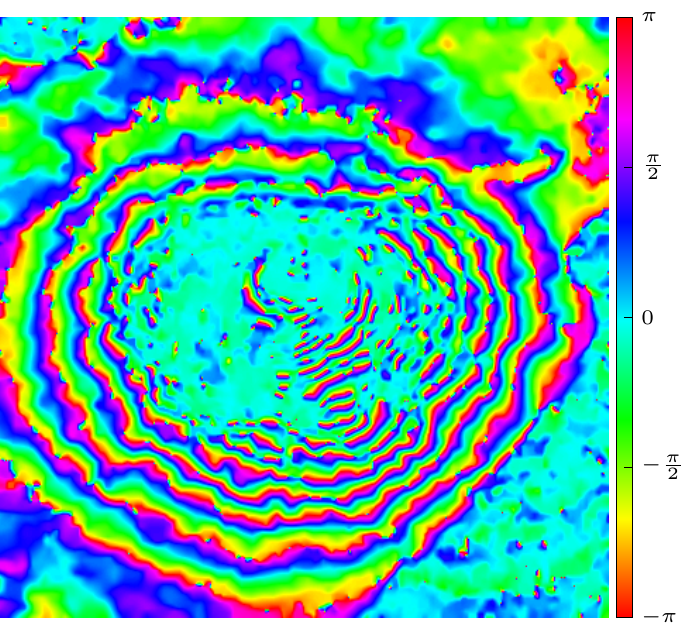}
			\caption{Inpainted and denoised version,\\where 20\% of data was lost.}\label{subfig:Vesuv:Lost}
		\end{subfigure}
		\caption{Real data of Mount Vesuvius. We compare a pure denoising approach in~\subref{subfig:Vesuv:Denoised} with a combined inpainting and denoising approach in~\subref{subfig:Vesuv:Lost}, where 20\% of the data was lost before the inpainting and denoising process.}\label{fig:vesuv}
\end{figure}

For both examples the computation takes about 43 seconds for first order
differences and 55 seconds for the combined cases, respectively. We used
\(k=2000\) iterations as a stopping criterion for the CPPA from
Sec.~\ref{sec:cpp}.
%
%
\paragraph{Inpainting in the presence of noise.}
In real world measurements data are often noisy. If these data are also
partially lost, we employ the
model~\eqref{eq:2DTVfunctionalInpaintingWithNoise}. As an example we consider
the measurement of elevation using InSAR~\cite{BRF00,MF98}. In particular, we
consider phase valued measured data of Mount
Vesuvius~\cite{RPG97}\footnote{online available at~
\url{%
	https://earth.esa.int/workshops/ers97/program-details/speeches/rocca-et-al/} 
}.
We compare denoising with simultaneously inpainting and denoising. To this end,
we randomly destroyed 20\% of the data items. The results without and with lost
data are shown in Fig.~\ref{fig:vesuv}\,(\subref{subfig:Vesuv:Denoised})
and~(\subref{subfig:Vesuv:Lost}), respectively. For the inpainting version the
parameters used in Fig.~\ref{fig:vesuv}\,(\subref{subfig:Vesuv:Denoised}),
\(\alpha=\frac{1}{4}(1,1,1,1)^\tT\), \(\beta = \frac{3}{4}(1,1)^\tT\) and
\(\gamma=\frac{3}{4}\), were multiplied by \(2\). The combined approach of
simultaneously inpainting and denoising introduces a few more artifacts than
pure denoising; cf. the middle and top right area. However, both results are of
comparable quality in smooth regions, e.g., the plateau in the bottom left.
%
\section{Conclusions} \label{sec:conclusions}
\enlargethispage{\baselineskip}
We proposed an inpainting model for cyclic data which involves our
recently established second order cyclic differences. 
Since there are analytical expressions for the proximal mappings of these 
differences we suggested a CPP algorithm together with a strategy for choosing
the cycles to compute a minimizer of the corresponding functionals efficiently.
There is large room for improvements and future work.

We want to apply our second order cyclic differences to other image restoration
tasks such as, e.g., deblurring and investigate other couplings of first and
second order differences.
It is possible to generalize our geometrically driven definition of second
order differences to higher dimensional spheres and also to general manifolds.
We want to use such generalization for image processing tasks of general manifold-valued data.
%
%
\section*{Appendix}\label{sec:appendix}
The proximal mappings from Theorem~\ref{lem:proxy_b1} can be efficiently
applied in parallel if they act on distinct data. This reduces the cycle
length~\(c\) of the CPPA from Sec.~\ref{sec:cpp} tremendously and provides
an efficient, parallel implementation. Especially, the cycle length is
independent of the inpainting area \(\Omega\) or the image domain \(\Omega_0\)
and depends only on the number of dissimilar differences used. We present a
specific splitting in the CPPA of our inpainting
model~\eqref{eq:2DTVfunctionalInpainting} given by
\begin{equation} \label{eq:splittingJ}
	J(x) = \alpha \operatorname{TV}_{1}^{\Omega}(x)
		+ \beta \operatorname{TV}_{2}^{\Omega}(x)
		+ \gamma \operatorname{TV}_{1,1}^{\Omega}(x) 
\end{equation}
with the constraints $x_{i,j}=f_{i,j}$ on  $\bar \Omega$.
We write 
\begin{equation} \label{splitting_18}
J = \sum_{l=1}^{18} J_l
\end{equation}
with summands $J_l$ given by the subsequent explanation.
We start with the~$\alpha \operatorname{TV}_{1}^{\Omega}(x)$ term
and first consider the horizontal
summand~$\alpha_1 \sum_{(i,j)} d(x_{i,j},x_{i+1,j})$.
We split this sum into an even and an odd part $J_1$ and $J_2$, more precisely
\begin{equation*}
	\alpha_1 \sum_{(i,j)} d(x_{i,j},x_{i+1,j})  = J_1 + J_2,	
\end{equation*}
where
\begin{equation*}
	J_1+J_2 :=\alpha_1 \sum_{(i,j)} d(x_{2i,j},x_{2i+1,j}) +
	\alpha_1 \sum_{(i,j)} d(x_{2i+1,j},x_{2i+2,j}),
\end{equation*}
with the restriction to the summands as in Sec.~\ref{sec:model}. This means,
that for each item \((i,j)\) in the sum, the corresponding index of at least
one of the arguments~\(x_{2i,j}\), \(x_{2i+1,j}\) is in \(\Omega\).  For the
vertical as well as for the diagonal summands
in~$\alpha \operatorname{TV}_{1}^{\Omega}(x)$
we proceed analogously to obtain the splitting functionals $J_3,\ldots,J_8$.

Next, we consider the $\beta \operatorname{TV}_{2}^{\Omega}(x)$ term with its
first (horizontal) summand given
by~$\beta_1 \sum_{(i,j)} d_2(x_{i-1,j},x_{i,j},x_{i+1,j})$.
We decompose this summand into three sums $J_9,J_{10},J_{11}$ given by 
\begin{align*}
   J_9   &= \beta_1 \sum_{(i,j)} d_2(x_{3i-1,j},x_{3i,j},x_{3i+1,j}),\\
   J_{10}  &= \beta_1 \sum_{(i,j)} d_2(x_{3i,j},x_{3i+1,j},x_{3i+2,j}),\\
   J_{11}  &= \beta_1 \sum_{(i,j)} d_2(x_{3i+1,j},x_{3i+2,j},x_{3i+3,j}),
\end{align*}   
again, with the restriction to the summands as in Sec.~\ref{sec:model}.
For the vertical summand in $\beta \operatorname{TV}_{2}^{\Omega}(x)$ we
proceed analogously to obtain $J_{12},\ldots,J_{14}.$

It remains to split the term $\gamma \operatorname{TV}_{1,1}^{\Omega}(x)$ into
four functionals~$J_{15},\ldots,J_{18}$ as follows
\begin{align*}
   J_{15}   &= \gamma \sum_{(i,j)} d_{1,1}(x_{2i,2j},x_{2i+1,2j},x_{2i,2j+1},x_{2i+1,2j+1}),\\
   J_{16}   &= \gamma \sum_{(i,j)} d_{1,1}(x_{2i+1,2j},x_{2i+2,2j},x_{2i+1,2j+1},x_{2i+2,2j+1}),\\
   J_{17}   &= \gamma \sum_{(i,j)} d_{1,1}(x_{2i,2j+1},x_{2i+1,2j+1},x_{2i,2j+2},x_{2i+1,2j+2}),\\
   J_{18}   &= \gamma \sum_{(i,j)} d_{1,1}(x_{2i+1,2j+1},x_{2i+2,2j+1},x_{2i+1,2j+2},x_{2i+2,2j+2}).\\
\end{align*}
Each summation is again restricted to those terms where at least one index of
an argument of \(d_{1,1}\) is in \(\Omega\).
For $J_1,\ldots,J_{18}$ the corresponding proximal mapping can be explicitly
computed and the cycle length \(c=18\) is independent of the cardinality
of~\(\Omega\) or \(\Omega_0\).
After application of the proximal mapping of each $J_i$ we set
$x_{i,j}=f_{i,j}$ on  $\bar \Omega$ to fulfill the respective constraint, which
is the same as performing a projection.
%

\begin{thebibliography}{10}

\bibitem{AF13}
M.~Almeida and M.~Figueiredo.
\newblock Deconvolving images with unknown boundaries using the alternating
  direction method of multipliers.
\newblock {\em IEEE Trans. on Image Process.}, 22(8):3074--3086, 2013.

\bibitem{Bac13}
M.~Ba{\v c}{\'a}k.
\newblock {The proximal point algorithm in metric spaces}.
\newblock {\em Isr. J. Math.}, 194(2):689--701, 2013.

\bibitem{Bac13a}
M.~Ba{\v c}{\'a}k.
\newblock Computing medians and means in {H}adamard spaces.
\newblock {\em SIAM J. Optim.}, 2014.
\newblock to appear.

\bibitem{BBCSV01}
C.~Ballester, M.~Bertalmio, V.~Caselles, G.~Sapiro, and J.~Verdera.
\newblock Filling in by joint interpolation of vector fields and gray levels.
\newblock {\em IEEE Trans. Image Process.}, 10(8):1200--1211, 2001.

\bibitem{BLSW14}
R.~Bergmann, F.~Laus, G.~Steidl, and A.~Weinmann.
\newblock Second order differences of cyclic data and applications in
  variational denoising.
\newblock {\em Preprint}, 2014.

\bibitem{BSCB00}
M.~Bertalmío, G.~Sapiro, V.~Caselles, and C.~Ballester.
\newblock Image inpainting.
\newblock In {\em Proceedings of SIGGRAPH 2000}, pages 417--424, New Orleans,
  USA, 2000.

\bibitem{Ber10}
D.~P. Bertsekas.
\newblock Incremental gradient, subgradient, and proximal methods for convex
  optimization: a survey.
\newblock Technical Report LIDS-P-2848, Laboratory for Information and Decision
  Systems, MIT, Cambridge, MA, 2010.

\bibitem{Ber11}
D.~P. Bertsekas.
\newblock {Incremental proximal methods for large scale convex optimization}.
\newblock {\em Math. Program., Ser. B}, 129(2):163--195, 2011.

\bibitem{BKP09}
K.~Bredies, K.~Kunisch, and T.~Pock.
\newblock {Total generalized variation}.
\newblock {\em SIAM J. Imaging Sci.}, 3(3):1--42, 2009.

\bibitem{BBCS10}
A.~Bugeau, M.~Bertalmío, V.~Caselles, and G.~Sapiro.
\newblock A comprehensive framework for image inpainting.
\newblock {\em IEEE Trans. Signal Process.}, 19:2634--2645, 2010.

\bibitem{BRF00}
R.~B{\"u}rgmann, P.~A. Rosen, and E.~J. Fielding.
\newblock Synthetic aperture radar interferometry to measure
  earth{\textquoteright}s surface topography and its deformation.
\newblock {\em Annu. Rev. Earth Planet. Sci.}, 28(1):169--209, 2000.

\bibitem{CDOS12}
J.-F. Cai, B.~Dong, S.~Osher, and Z.~Shen.
\newblock Image restoration: Total variation, wavelet frames, and beyond.
\newblock {\em J. Amer. Math. Soc.}, 25(4):1033--1089, 2012.

\bibitem{caselles1998axiomatic}
V.~Caselles, J.-M. Morel, and C.~Sbert.
\newblock An axiomatic approach to image interpolation.
\newblock {\em IEEE Trans. on Image Process.}, 7(3):376--386, 1998.

\bibitem{CL97}
A.~Chambolle and P.-L. Lions.
\newblock {Image recovery via total variation minimization and related
  problems}.
\newblock {\em Numer. Math.}, 76(2):167--188, 1997.

\bibitem{CS01}
T.~Chan and J.~Shen.
\newblock Local inpainting models and {TV} inpainting.
\newblock {\em SIAM J. Appl. Math.}, 62(3):1019--1043, 2001.

\bibitem{CMM00}
T.~F. Chan, A.~Marquina, and P.~Mulet.
\newblock {High-order total variation-based image restoration}.
\newblock {\em SIAM J. Sci. Comput.}, 22(2):503--516, 2000.

\bibitem{chan2005image}
T.~F. Chan and J.~Shen.
\newblock {\em Image Processing and Analysis: Variational, PDE, Wavelet, and
  Stochastic Methods}.
\newblock SIAM, 2005.

\bibitem{CTDF04}
C.~Chefd'Hotel, D.~Tschumperl{\'e}, R.~Deriche, and O.~Faugeras.
\newblock {Regularizing flows for constrained matrix-valued images}.
\newblock {\em J. Math. Imaging Vis.}, 20(1-2):147--162, 2004.

\bibitem{DWB09}
S.~Didas, J.~Weickert, and B.~Burgeth.
\newblock {Properties of higher order nonlinear diffusion filtering}.
\newblock {\em J. Math. Imaging Vis.}, 35:208--226, 2009.

\bibitem{FO02}
O.~P. Ferreira and P.~R. Oliveira.
\newblock {Proximal point algorithm on {R}iemannian manifolds}.
\newblock {\em Optimization}, 51(2):257--270, 2002.

\bibitem{fisher95}
N.~I. Fisher.
\newblock {Statistical Analysis of Circular Data}.
\newblock Cambridge University Press, 1995.

\bibitem{Fle13}
P.~Fletcher.
\newblock Geodesic regression and the theory of least squares on {R}iemannian
  manifolds.
\newblock {\em Int. J. Comput. Vision}, 105(2):171--185, 2013.

\bibitem{FJ07}
P.~Fletcher and S.~Joshi.
\newblock {Riemannian geometry for the statistical analysis of diffusion tensor
  data}.
\newblock {\em Signal Process.}, 87(2):250--262, 2007.

\bibitem{GP98}
D.~C. Ghiglia and M.~D. Pritt.
\newblock {\em {Two-dimensional phase unwrapping: theory, algorithms, and
  software}}.
\newblock Wiley, 1998.

\bibitem{GMS93}
M.~Giaquinta, G.~Modica, and J.~Sou{\v c}ek.
\newblock Variational problems for maps of bounded variation with values in
  {$S^1$}.
\newblock {\em Calc. Var.}, 1(1):87--121, 1993.

\bibitem{GM06}
M.~Giaquinta and D.~Mucci.
\newblock The {BV}-energy of maps into a manifold: relaxation and density
  results.
\newblock {\em Ann. Sc. Norm. Super. Pisa Cl. Sci.}, 5(4):483--548, 2006.

\bibitem{GM07}
M.~Giaquinta and D.~Mucci.
\newblock {Maps of bounded variation with values into a manifold: total
  variation and relaxed energy}.
\newblock {\em Pure Appl. Math. Q.}, 3(2):513--538, 2007.

\bibitem{GHS13}
P.~Grohs, H.~Hardering, and O.~Sander.
\newblock Optimal a priori discretization error bounds for geodesic finite
  elements.
\newblock Technical Report 2013-16, Seminar for Applied Mathematics, ETH
  Z{\"u}rich, Switzerland, 2013.

\bibitem{GW09}
P.~Grohs and J.~Wallner.
\newblock {Interpolatory wavelets for manifold-valued data}.
\newblock {\em Appl. Comput. Harmon. Anal.}, 27(3):325--333, 2009.

\bibitem{guillemot2014image}
C.~Guillemot and O.~Le~Meur.
\newblock Image inpainting: Overview and recent advances.
\newblock {\em IEEE Signal Process. Mag.}, 31(1):127--144, 2014.

\bibitem{HS06}
W.~Hinterberger and O.~Scherzer.
\newblock Variational methods on the space of functions of bounded {H}essian
  for convexification and denoising.
\newblock {\em Computing}, 76(1):109--133, 2006.

\bibitem{JS2001}
S.~R. Jammalamadaka and A.~SenGupta.
\newblock {\em {Topics in Circular Statistics}}.
\newblock World Scientific Publishing Company, 2001.

\bibitem{LBU2012}
S.~Lefkimmiatis, A.~Bourquard, and M.~Unser.
\newblock Hessian-based norm regularization for image restoration with
  biomedical applications.
\newblock {\em IEEE Trans. Image Process.}, 21(3):983--995, 2012.

\bibitem{LSKC13}
J.~Lellmann, E.~Strekalovskiy, S.~Koetter, and D.~Cremers.
\newblock {Total variation regularization for functions with values in a
  manifold}.
\newblock In {\em IEEE ICCV 2013}, pages 2944--2951, 2013.

\bibitem{LLT03}
M.~Lysaker, A.~Lundervold, and X.-C. Tai.
\newblock {Noise removal using fourth-order partial differential equations with
  applications to medical magnetic resonance images in space and time}.
\newblock {\em IEEE Trans. Image Process.}, 12(12):1579--1590, 2003.

\bibitem{marz2011image}
T.~M{\"a}rz.
\newblock Image inpainting based on coherence transport with adapted distance
  functions.
\newblock {\em SIAM J. Imaging Sci.}, 4(4):981--1000, 2011.

\bibitem{marz2013well}
T.~M{\"a}rz.
\newblock A well-posedness framework for inpainting based on coherence
  transport.
\newblock {\em Found. Comput. Math.}, 2014.
\newblock to appear.

\bibitem{masnou1998level}
S.~Masnou and J.-M. Morel.
\newblock Level lines based disocclusion.
\newblock In {\em IEEE ICIP 1998}, pages 259--263, 1998.

\bibitem{MF98}
D.~Massonnet and K.~L. Feigl.
\newblock {R}adar interferometry and its application to changes in the
  {E}arth's surface.
\newblock {\em Rev. Geophys.}, 36(4):441--500, 1998.

\bibitem{PS13}
K.~Papafitsoros and C.~B. Sch{\"o}nlieb.
\newblock {A combined first and second order variational approach for image
  reconstruction}.
\newblock {\em J. Math. Imaging Vis.}, 2(48):308--338, 2014.

\bibitem{Pen06}
X.~Pennec.
\newblock Intrinsic statistics on {R}iemannian manifolds: Basic tools for
  geometric measurements.
\newblock {\em J. Math. Imaging Vis.}, 25(1):127--154, 2006.

\bibitem{RDSDS05}
I.~U. Rahman, I.~Drori, V.~C. Stodden, and D.~L. Donoho.
\newblock {Multiscale representations for manifold-valued data}.
\newblock {\em Multiscale Model. Simul.}, 4(4):1201--1232, 2005.

\bibitem{RPG97}
F.~Rocca, C.~Prati, and A.~M. Guarnieri.
\newblock Possibilities and limits of {SAR} interferometry.
\newblock In {\em Proc. Int. Conf. Image Process. Techn. 1996}, pages 15--26,
  1997.

\bibitem{Roc76}
R.~T. Rockafellar.
\newblock Monotone operators and the proximal point algorithm.
\newblock {\em SIAM J. Control Optim.}, 14(5):877--898, 1976.

\bibitem{ROF92}
L.~I. Rudin, S.~Osher, and E.~Fatemi.
\newblock {Nonlinear total variation based noise removal algorithms}.
\newblock {\em Physica D.}, 60(1):259--268, 1992.

\bibitem{Sche98}
O.~Scherzer.
\newblock {Denoising with higher order derivatives of bounded variation and an
  application to parameter estimation}.
\newblock {\em Computing}, 60:1--27, 1998.

\bibitem{SS08}
S.~Setzer and G.~Steidl.
\newblock Variational methods with higher order derivatives in image
  processing.
\newblock In {\em Approximation Theory XII: San Antonio 2007}, pages 360--385,
  2008.

\bibitem{SST11}
S.~Setzer, G.~Steidl, and T.~Teuber.
\newblock {Infimal convolution regularizations with discrete l1-type
  functionals}.
\newblock {\em Commun. Math. Sci.}, 9(3):797--872, 2011.

\bibitem{SC11}
E.~Strekalovskiy and D.~Cremers.
\newblock {Total variation for cyclic structures: Convex relaxation and
  efficient minimization}.
\newblock In {\em IEEE CVPR 2011}, pages 1905--1911, 2011.

\bibitem{CS13}
E.~Strekalovskiy and D.~Cremers.
\newblock Total cyclic variation and generalizations.
\newblock {\em J. Math. Imaging Vis.}, 47(3):258--277, 2013.

\bibitem{VBK13}
T.~Valkonen, K.~Bredies, and F.~Knoll.
\newblock Total generalized variation in diffusion tensor imaging.
\newblock {\em SIAM J. Imag. Sci.}, 6(1):487--525, 2013.

\bibitem{Wein12}
A.~Weinmann.
\newblock {Interpolatory multiscale representation for functions between
  manifolds}.
\newblock {\em SIAM J. Math. Anal.}, 44(1):162--191, 2012.

\bibitem{WDS2013}
A.~Weinmann, L.~Demaret, and M.~Storath.
\newblock Total variation regularization for manifold-valued data.
\newblock {\em Preprint}, 2013.

\end{thebibliography}
%

%
\end{document}